\documentclass[aos,MSNbibl,seceqn,dvips]{arximspdf}
\usepackage{dcolumn}

\doi{10.1214/13-AOS1090} 
\volume{41}
\issue{2}
\pubyear{2013}
\firstpage{484}
\lastpage{507}

\makeatletter

\newcolumntype{d}[1]{D{.}{.}{#1}}

\newcommand{\cal}{\mathcal}
\newcommand{\E}{\mathrm{E}}
\newcommand{\dto}{\stackrel{d}{\rightarrow}}
\renewcommand{\P}{\mathrm{P}}
\newcommand{\I}{\mathrm{I}}
\newcommand{\Binomial}{\operatorname{Binomial}}
\def\Larrow{\stackrel{L}{\rightarrow}}
\def\Parrow{\stackrel{P}{\rightarrow}}

\newtheorem{theorem}{Theorem}[section]
\newtheorem{lemma}{Lemma}[section]

\newproclaim{example}{Example}[section]
\newproclaim{remark}{Remark}[section]
\newproclaim{examplecontd}{Example}[section]

\makeatother

\begin{document}
\begin{frontmatter}

\title{Exact and asymptotically robust permutation tests\thanksref{T1}}
\runtitle{Exact and asymptotically robust permutation tests}

\thankstext{T1}{Supported by NSF Grant DMS-07-07085.}

\begin{aug}
\author[A]{\fnms{EunYi} \snm{Chung}\ead[label=e1]{eunyi@stanford.edu}}
\and
\author[B]{\fnms{Joseph P.} \snm{Romano}\corref{}\ead[label=e2]{romano@stanford.edu}}
\runauthor{E. Chung and J. P. Romano}
\affiliation{Stanford University}
\address[A]{Department of Economics\\
Stanford University\\
Stanford, California 94305-6072\\
USA\\
\printead{e1}} 
\address[B]{Departments of Statistics and Economics\\
Stanford University\\
Stanford, California 94305-4068\\
USA\\
\printead{e2}}
\end{aug}

\received{\smonth{7} \syear{2012}}
\revised{\smonth{12} \syear{2012}}

%
\begin{abstract}
Given independent samples from $P$ and $Q$, two-sample permutation
tests allow one to construct exact level tests when the null hypothesis
is $P = Q$. On the other hand, when comparing or testing particular
parameters $\theta$ of $P$ and $Q$, such as their means or medians,
permutation tests need not be level~$\alpha$, or even approximately
level $\alpha$ in large samples. Under very weak assumptions for
comparing estimators, we provide a general test procedure whereby the
asymptotic validity of the permutation test holds while retaining the
\textit{exact} rejection probability $\alpha$ in finite samples when the
underlying distributions are identical. The ideas are broadly
applicable and special attention is given to the $k$-sample problem of
comparing general parameters, whereby a permutation test is constructed
which is exact level $\alpha$ under the hypothesis of identical
distributions, but has asymptotic rejection probability $\alpha$ under
the more general null hypothesis of equality of parameters. A Monte
Carlo simulation study is performed as well. A quite general theory is
possible based on a coupling construction, as well as a key contiguity
argument for the multinomial and multivariate hypergeometric
distributions.
\end{abstract}

%
\begin{keyword}[class=AMS]
\kwd[Primary ]{62E20}
\kwd[; secondary ]{62G10}
\end{keyword}
\begin{keyword}
\kwd{Behrens--Fisher problem}
\kwd{coupling}
\kwd{permutation test}
\end{keyword}

\end{frontmatter}

\section{Introduction}\label{secintroduction}

In this article, we consider the behavior of two-sample (and later also
$k$-sample) permutation tests
for testing problems when the fundamental assumption of identical
distributions need not hold. Assume $X_1,\ldots, X_m$ are i.i.d.
according to a probability distribution $P$, and independently, $Y_1,\ldots, Y_n$
are i.i.d.~$Q$. The underlying model specifies a family of pairs of
distributions
$(P,Q)$ in some space $\Omega$. For the problems considered here,
$\Omega$
specifies a nonparametric model, such as the set of all pairs of distributions.
Let $N = m+n$, and write
%
\begin{equation}
\label{equationZN} Z = (Z_1,\ldots, Z_N ) = (
X_1,\ldots, X_m, Y_1,\ldots,
Y_n ).
\end{equation}
Let $\bar\Omega= \{ (P,Q)\dvtx  P = Q \}$. Under the assumption $(P,Q)
\in\bar\Omega$, the joint distribution of $(Z_1,\ldots, Z_N)$ is the
same as $(Z_{\pi(1)},\ldots, Z_{\pi(N) } )$, where $(
\pi(1),\ldots,\break
\pi(N) )$ is any permutation of $\{ 1,\ldots, N \}$. It follows that,
when testing any null hypothesis $H_0\dvtx  (P, Q) \in\Omega_0$, where
$\Omega_0 \subset\bar\Omega$, then an exact level $\alpha$ test can be
constructed by a permutation test. To review how, let $\mathbf{G}_N$
denote the set of all permutations $\pi$ of $\{ 1,\ldots, N \}$. Then,
given any test statistic $T_{m,n} = T_{m,n} ( Z_1,\ldots, Z_N )$,
recompute $T_{m,n}$ for all permutations $\pi$; that is, compute
$T_{m,n} ( Z_{\pi(1)},\ldots, Z_{\pi(N)} )$ for all
$\pi\in\mathbf{G}_N$, and let their ordered values be\looseness=-1
\[
T_{m,n}^{(1)} \le T_{m,n}^{(2)} \le\cdots\le
T_{m,n}^{( N! )}.
\]\looseness=0
Fix a nominal level $\alpha$, $0 < \alpha< 1$, and let $k$ be defined
by $k = N! - \lbrack\alpha N! \rbrack$,
where $\lbrack\alpha N! \rbrack$ denotes the largest integer less
than or equal to
$ \alpha N!$.
Let $M^+ (z)$ and $M^0 (z)$ be the number of values $T_{m,n}^{(j)} (z)$
$( j=1,\ldots, N!)$ which are greater than $T^{(k)} (z)$
and equal to $T^{(k)} (z)$, respectively.
Set
\[
a(z) = \frac{ \alpha N! - M^+ (z) }{ M^0 (z) }.
\]

Define the randomization test function $\phi(Z)$ to be equal to
1, $a(Z)$ or 0 according to whether $T_{m,n} (Z) > T_{m,n}^{(k)} (Z)$,
$T_{m,n} (X) = T^{(k)} (Z)$ or $T_{m,n} (Z) < T^{(k)} (Z)$, respectively.
Then, under any $(P,Q) \in\bar\Omega$,
\[
E_{P,Q} \bigl[ \phi(X_1,\ldots, X_m,
Y_1,\ldots, Y_n ) \bigr] = \alpha.
\]

Also, define the permutation distribution as
%
\begin{equation}
\label{equationnewperm} \hat R_{m,n}^T (t) =
\frac{1}{N!} \sum_{\pi\in\mathbf{G}_N} I \bigl\{
T_{m,n} ( Z_{\pi(1)},\ldots, Z_{\pi(N)} ) \le t \bigr\}.
\end{equation}
Roughly speaking (after accounting for discreteness), the permutation
test rejects $H_0$ if the test statistic $T_{m,n}$ exceeds
$T_{m,n}^{(k)}$, or a $1- \alpha$ quantile of this permutation
distribution.

It may be helpful to consider an alternative description of the
permutation distribution given in
(\ref{equationnewperm}). As a shorthand, for any $\pi\in\mathbf{G}_N$, let
$Z_{\pi} = (Z_{ \pi(1)},\ldots, Z_{\pi(N)} )$.
Let $\Pi$ denote a random permutation, uniformly distributed
over $\mathbf{G}_N$. Then, $T_{m,n} ( Z_{\Pi})$ denotes the random
variable that
evaluates the test statistic, not at the original data $Z$, but at a
randomly permuted
data set $Z_{\Pi}$. The permutation distribution $\hat R_{m,n}^T (
\cdot)$ given
in (\ref{equationnewperm}) is evidently the conditional distribution
of $T_{m,n} (Z_{\Pi} ) $
given $Z$, because conditional on the data $Z$, $T_{m,n} ( Z_{\Pi} )$
is equally likely
to be any of $T_{m,n} (Z_{\pi} )$ among $\pi\in\mathbf{G}_N$.
The asymptotic behavior of this (conditional) distribution $\hat
R_{m,n}^T ( \cdot)$ is the key
to establishing properties of the permutation test.

Although the rejection probability of the permutation test is exactly
$\alpha$ when \mbox{$P = Q$}, problems arise
if $\Omega_0$ is strictly bigger than $\bar\Omega$. Since a
transformed permuted
data set no longer has the same distribution as the original data set,
the argument leading
to the construction of an $\alpha$ level test fails, and faulty
inferences can occur.

To be concrete, if we are interested in testing equality of means, for
example, then $\Omega_0 = \{(P,Q)\dvtx \mu(P) = \mu(Q) \}$ which, of
course, is strictly bigger than $\bar\Omega$. So, consider
constructing a permutation test based on
the difference of sample means
%
\begin{equation}
\label{equationtmn} T_{m,n} = \sqrt{N} ( \bar X_m - \bar
Y_n ).
\end{equation}
Note that we are not taking the absolute difference, so that the test
is one-sided, as we are rejecting for large positive values of the difference.
First of all, we are not concerned about testing $\bar\Omega= \{
(P,Q)\dvtx P = Q \}$, but something bigger than $\bar\Omega$. However, we
underscore the point that a test statistic (\ref{equationtmn}) is not
appropriate for testing $\bar\Omega$ without further assumptions
because the test clearly will not have any power
against distributions $P$ and $Q$ whose means are identical but $P \ne Q$.

The permutation test
based on the difference of sample means is only appropriate
as a test of equality of population means. However, the permutation test
no longer controls the level of the test, even in large samples.
As is well known (Romano \cite{r27}), the permutation test possesses a
certain asymptotic
robustness as a test of difference in means if $m/n \to1$ as $n \to
\infty$, or the underlying variances of $P$ and $Q$ are equal, in
the sense that the rejection probability under the null hypothesis
of equal means tends to the nominal level. Without equal variances or
comparable sample sizes, the
rejection probability can be much larger than the nominal level,
which is a concern. Because of
the lack of robustness and the increased probability of a type 1
error, rejection of the null may incorrectly be interpreted as
rejection of equal means, when in fact it is caused by unequal
variances and unequal sample sizes. Even more alarming is the
possibility of rejecting a two-sided null hypothesis when observing a
positive large
difference with the accompanying inference
that mean difference is positive when in fact the difference in means is
negative, a~type 3 error or directional error.
Indeed, if for some $P$ and $Q$ with equal means the rejection
probability is, say, $\gamma\gg\alpha$, then it follows by continuity
that the rejection probability under some $P$ and $Q$
with negative mean difference will be nearly $\gamma$ as well, where
one would conclude
that the mean difference is actually positive.
Further note that there is also the possibility that the
rejection probability can be much less than the nominal level, which
by continuity implies the test is biased and has little power of
detecting a true difference in means, or large type 2 error.

The situation is even worse when basing a test on a difference in
sample medians,
in the sense that regardless of sample sizes, the asymptotic rejection
probability
of the permutation test will be $\alpha$ under very stringent
conditions, which
essentially means only in the case where the underlying distributions
are the same.

However, in a very insightful paper in the context of random censoring models,
Neuhaus \cite{r21} realized that by proper studentization of a test
statistic, the permutation test can result in asymptotically valid
inference even when
the underlying distributions are not the same.
This result has been extended to other \mbox{specific} problems, such as
comparing means by Janssen \cite{r9} and certain linear statistics in
Janssen \cite{r10} (including the Wilcoxon statistic without ties),
variances by Pauly \cite{r23} and the two-sample Wilcoxon test by
Neubert and Brunner \cite{r20} (where ties are allowed).
Other results on permutation tests are presented in Janssen \cite
{r11}, Janssen and Pauls \cite{r12}, Janssen and Pauls \cite{r13} and
Janssen and Pauly \cite{r14}. The recent paper by
Omelka and Pauly \cite{r22} compares correlations by permutation
tests, which is a special
case of our general results.
Note that the importance of studentization when bootstrapping is well
known; see
Hall and Wilson \cite{r7} and Delaigle et al. \cite{r3} (though its
role for bootstrap is to obtain higher order accuracy
while in the context here first order accuracy can fail without studentization).

The goal of this paper is to obtain a quite general result of the same
phenomenon.
That is, when basing a permutation test using some test statistic as a
test of a
parameter (usually a difference of parameters associated with marginal
distributions),
we would like to retain the exactness property when $P =Q$, and also have
the asymptotic rejection probability be $\alpha$ for the more general
null hypothesis
specifying the parameter (such as the difference being zero).
Of course, there are many alternatives to getting asymptotic tests,
such as the bootstrap or subsampling. However, we do not wish
to give up the exactness property under $P = Q$, and resampling methods
do not have such finite sample guarantees.
The main problem becomes: what is the asymptotic behavior of $\hat
R_{m,n}^T ( \cdot)$
defined in (\ref{equationnewperm}) for general test statistic
sequences $T_{m,n}$
when the underlying distributions differ.
Only for suitable test statistics is it possible to achieve both finite sample
exactness when the underlying distributions are equal, but also
maintain a large sample rejection probability near the nominal level
when the underlying distributions need not be equal.
In this sense, our results are both exact and asymptotically robust
for heterogeneous populations.

This paper provides a framework for testing a parameter that depends on
$P$ and~$Q$ (and later on $k$ underlying distributions $P_i$ for $i =
1,\ldots, k$). We construct a general test procedure where the asymptotic
validity of the permutation test holds in a general setting.
Assuming that estimators are asymptotically linear and consistent estimators
are available for their asymptotic variance, we provide a test
that has asymptotic rejection probability equal to the nominal level
$\alpha$,
but still retains the exact rejection probability of $\alpha$ in
finite samples if $P = Q$ in Section~\ref{secstudentized}.
It is not even required that the estimators are based on differentiable
functionals,
and some methods like the bootstrap would not necessarily be even
asymptotically valid
under such conditions, let alone retain the finite sample exactness property
when $P=Q$. In Section~\ref{sectiongeneralization}, generalizations
of the results are discussed with a special attention to the more
general $k$-sample problem of comparing general parameters.
Furthermore, Monte Carlo simulation studies illustrating our results
are presented in Section~\ref{secsimulation}. The arguments of the
paper are quite different from Janssen and previous authors, and hold
under great generality. For example, they immediately apply to
comparing means,
variances or medians. The key idea is to show that the permutation distribution
behaves like the unconditional distribution of the test statistic
when all $N$ observations are i.i.d. from the mixture distribution $pP
+ (1- p)Q$,
where $p$ is such that $m/N \to p$. This seems intuitive because the
permutation distribution
permutes the observations so that a permuted sample is almost like a sample
from the mixture distribution. In order to make this idea precise, a coupling
argument is given in Section~\ref{sectioncoupling}. Of course, the
permutation distribution depends on
all permuted samples (for a given original data set). But even for one permuted
data set, it cannot exactly be viewed as a sample from $pP + (1-p)Q$.
Indeed, the first $m$ observations from the mixture would include $B_m$
observations from $P$ and the rest from $Q$, where $B_m$ has the binomial
distribution based on $m$ trials and success probability $p$.
On the other hand, for a permuted sample, if $H_m$ denotes the number of
observations from $P$, then $H_m$ has the hypergeometric distribution
with mean $mp$. The key argument that allows for such a general
result concerns the contiguity of the distributions of $B_m$ and $H_m$.
Section~\ref{sectioningredients} highlights the main technical ideas
required for the proofs. All proofs are deferred to the supplementary
appendix \cite{r2}.

\section{Robust studentized two-sample test}\label{secstudentized}
In this section, we consider the general problem of inference from the
permutation
distribution
when comparing parameters from two populations.
Specifically, assume $X_1,\ldots, X_m$ are i.i.d. $P$ and, independently,
$Y_1,\ldots, Y_n$ are i.i.d. $Q$.
Let $\theta( \cdot)$ be a real-valued parameter, defined on some space
of distributions $\cal P$.
The problem is to test the null hypothesis
%
\begin{equation}
\label{equationnull} H_0\dvtx  \theta(P) = \theta(Q).
\end{equation}
Of course, when $P = Q$, one can construct permutation tests with exact
level $\alpha$. Unfortunately, if $P \ne Q$, the test need not be valid
in the sense that the probability of a type 1 error need not be $\alpha$
even asymptotically. Thus, our goal is to construct a procedure
that has asymptotic rejection probability equal to $\alpha$ quite generally,
but also retains the exactness property in finite samples when $P=Q$.

We will\vspace*{1pt} assume that estimators are available that are asymptotically linear.
Specifically, assume that, under $P$, there exists an estimator
$\hat\theta_m = \hat\theta_m ( X_1,\ldots, X_m )$ which satisfies
%
\begin{equation}
\label{equationxlinear} m^{1/2} \bigl[ \hat\theta_m -
\theta(P) \bigr] = \frac{1}{\sqrt{m}} \sum_{i=1}^m
f_P ( X_i ) + o_P (1).
\end{equation}
Similarly, we assume that, based on the $Y_j$ (under $Q$),
%
\begin{equation}
\label{equationylinear} n^{1/2} \bigl[ \hat\theta_n -
\theta(Q) \bigr] = \frac{1}{\sqrt{n}} \sum_{j=1}^n
f_Q ( Y_j) + o_Q (1).
\end{equation}
The functions determining the linear approximation $f_P$ and $f_Q$ can
of course depend on the underlying distributions. Different forms of
differentiability guarantee such linear expansions in the special case
when $\hat\theta_m$ takes the form of an empirical estimate $\theta(
\hat P_m)$, where $\hat P_m$ is the empirical measure constructed from
$X_1,\ldots, X_m$, but we will not need to assume such stronger
conditions. We will argue that our assumptions of asymptotic linearity
already imply\vspace*{2pt} a result about the permutation distribution corresponding
to the statistic $N^{1/2} [ \hat\theta_m ( X_1,\ldots, X_m ) -
\hat\theta_n ( Y_1,\ldots, Y_n )]$, without having to impose any
differentiability assumptions. However, we will assume the expansion
(\ref{equationxlinear}) holds not just for i.i.d. samples under $P$,
and also under $Q$, but also when sampling i.i.d. observations from the
mixture distribution $\bar P = pP + qQ$. This is a weak assumption and
replaces having to study the permutation distribution based on
variables that are no longer independent nor identically distributed
with a simple assumption about the behavior under an i.i.d. sequence.
Indeed, we will argue that in all cases, the permutation distribution
behaves asymptotically like the unconditional limiting sampling
distribution of the studied statistic sequence when sampling i.i.d.
observations from $\bar P$.

In the next two theorems, the behavior of the permutation distribution
is obtained. Note that it is not assumed that the null hypothesis
$\theta(P) = \theta(Q)$ necessarily holds. Indeed, the asymptotic
behavior of the permutation test under $P$ and $Q$ is the same as when
all observations are from the mixture distribution $\bar P = pP +
(1-p)Q$, where $p = \lim\frac{m}{N}$. Proofs of all the results in
Section~\ref{secstudentized} are presented along with proofs of the
results in Section~\ref{sectioningredients} in the supplementary
appendix \cite{r2}.


\begin{theorem}\label{theorembasic}
Assume $X_1,\ldots, X_m$ are i.i.d. $P$ and, independently,
$Y_1,\ldots, Y_n$ are i.i.d. $Q$. Consider testing the null hypothesis
(\ref{equationnull}) based on a test statistic of the form
\[
T_{m,n} = N^{1/2} \bigl[ \hat\theta_m
(X_1,\ldots, X_m ) - \hat\theta_n
(Y_1,\ldots, Y_n ) \bigr],
\]
where the estimators satisfy (\ref{equationxlinear}) and (\ref
{equationylinear}).
Further assume $E_P f_P (X_i ) = 0$ and
\[
0 < \E_P f_P^2 (X_i ) \equiv
\sigma^2 (P) < \infty
\]
and the same with $P$ replaced by $Q$. Let $m \to\infty$, $n \to
\infty$, with
$N = m+n$, $p_m = m/N$, $q_m = n/N$ and $p_m \to p \in(0,1)$ with
%
\begin{equation}
\label{equationmnmn} p_m - p = O\bigl( N^{-1/2} \bigr).
\end{equation}
Assume the estimator sequence
also satisfies (\ref{equationxlinear}) with $P$ replaced by $\bar P =
pP + qQ$
with $\sigma^2 ( \bar P ) < \infty$.

Then the permutation distribution of $T_{m,n}$ given by
(\ref{equationnewperm}) satisfies
\[
\sup_t \bigl| \hat R_{m,n}^T (t) - \Phi
\bigl( t / \tau(\bar P)\bigr)\bigr| \Parrow0,
\]
where
%
\begin{equation}
\label{equationtau2} \tau^2 (\bar P) = \frac{1}{p(1-p)}
\sigma^2 ( \bar P ).
\end{equation}
\end{theorem}

\begin{remark}\label{remark11}
Under $H_0$ given by (\ref{equationnull}), the true unconditional
sampling distribution of $T_{m,n}$
is asymptotically normal with mean 0 and variance
%
\begin{equation}
\label{equationtruevar} \frac{1}{p}\sigma^2 (P) +
\frac{1}{1-p} \sigma^2 (Q),
\end{equation}
which does not equal $\tau^2 ( \bar P )$ defined by (\ref
{equationtau2}) in general.
\end{remark}

\begin{example}[(Difference of means)]\label{exmeans}
As is well known, even for the case of comparing population means
by sample means, under the null hypothesis that $\theta(P) = \theta
(Q)$, equality of (\ref{equationtau2}) and (\ref{equationtruevar})
holds if and only if $p = 1/2$ or
$\sigma^2 (P) = \sigma^2 (Q)$.
\end{example}

\begin{example}[(Difference of medians)]\label{exmedians}
Let $F$ and $G$ denote the c.d.f.s corresponding to $P$ and $Q$.
Let $\theta(F)$ denote the median of $F$, that is, $\theta(F) = \inf
\{x\dvtx  F(x) \geq\frac{1}{2} \}$. Then it is well known (Serfling \cite
{r29}) that if $F$ is continuously differentiable at $\theta(P)$ with
derivative $F'$ (and the same with $F$ replaced by $G$), then
\[
m^{1/2}\bigl[\theta(\hat{P}_m) - \theta(P)\bigr] =
\frac{1}{\sqrt{m}} \sum_{i=1}^m
\frac{{1/2}- \I\{X_i \leq\theta(P)\}}{F' (\theta
(P))} + o_P(1)
\]
and similarly, with $P$ and $F$ replaced by $Q$ and $G$.
Thus, we can apply Theorem~\ref{theorembasic} and conclude that, when
$\theta(P) = \theta(Q) = \theta$, the permutation distribution of
$T_{m,n}$ is approximately a normal distribution with mean 0 and
variance
\[
\frac{1}{4p(1-p)[p F'(\theta) + (1-p) G'(\theta)]^2}
\]
in large samples. On the other hand, the true sampling distribution is
approximately a normal distribution with mean 0 and variance
%
\begin{equation}
\label{equationvvar} v^2 (P,Q) \equiv\frac{1}{p}
\frac{1}{4 [F'(\theta)]^2} + \frac
{1}{1-p}\frac{1}{4 [G'(\theta)]^2}.
\end{equation}
Thus the permutation distribution and the true unconditional sampling
distribution behave differently asymptotically unless $F' (\theta) =
G' (\theta)$ is satisfied. Since we do not assume $P = Q$, this
condition is a strong assumption. Hence, the permutation test for
testing equality of medians is generally not valid in the sense that
the rejection probability tends to a value that is far from the nominal
level $\alpha$.
\end{example}

The main goal now is to show how studentizing the test statistic
leads to a general correction.

\begin{theorem}\label{theoremgeneral} Assume the setup and conditions
of Theorem~\ref{theorembasic}. Further assume that $\hat\sigma_m (
X_1,\ldots, X_m )$ is a consistent estimator of $\sigma(P)$ when
$X_1,\ldots, X_m$
are i.i.d. $P$. Assume consistency also under $Q$ and $\bar P$, so
that $\hat\sigma_n (V_1,\ldots, V_n ) \Parrow\sigma( \bar P )$
as $n \to\infty$ when
the $V_i$ are i.i.d. $\bar P$.
Define the studentized test statistic
%
\begin{equation}
\label{equationstudentized} S_{m,n} = \frac{T_{m,n}}{V_{m,n}},
\end{equation}
where
\[
V_{m,n} = \sqrt{ \frac{N}{m}\hat\sigma_m^2
( X_1,\ldots, X_m ) + \frac{N}{n} \hat
\sigma_n^2 ( Y_1,\ldots, Y_n
)}
\]
and consider the permutation distribution defined in (\ref{equationnewperm})
with $T$ replaced by $S$.
Then
%
\begin{equation}
\label{equationratio} \sup_t \bigl| \hat R^S_{m,n}
(t) - \Phi(t) \bigr| \Parrow0.
\end{equation}
\end{theorem}

Thus the permutation distribution is asymptotically standard normal,
as is the true unconditional limiting distribution of the test
statistics $S_{m,n}$. Indeed, as mentioned in Remark~\ref{remark11},
the true unconditional limiting distribution of $T_{m,n}$ is
normal with mean 0 and variance given by (\ref{equationtruevar}).
But, when sampling $m$ observations from $P$ and $n$ from $Q$,
$V_{m,n}^2$ tends in probability to (\ref{equationtruevar}), and hence
the limiting distribution of $T_{m,n}$ is standard normal, the same
as that of the permutation distribution.

\begin{remark}\label{remarkpower}
As previously noted, Theorems~\ref{theorembasic} and
\ref{theoremgeneral} are true even if $\theta(P) \neq\theta(Q)$. If
$\theta(P) = \theta(Q)$, then the true sampling distribution of
$S_{m,n}$ and the permutation test become approximately the same.
However, if $\theta(P) \neq\theta(Q)$, then we get the power tending
to 1. Indeed, the critical value from the permutation distribution
asymptotically tends to a finite value $z_{1-\alpha}$ in probability,
while the test statistic tends to infinity in probability. Also, see
Remark~\ref{remarklocalpower} for local power.
\end{remark}
%
\begin{examplecontd}[(Continued)]
As proved by Janssen \cite{r9}, even when the underlying distributions
may have different variances and different sample sizes, permutation
tests based on studentized statistics
\[
S_{m,n} = \frac{N^{1/2}(\bar{X}_m - \bar{Y}_n)}{\sqrt{
{N}S^2_X/{m} + {N}S^2_Y/{n}}},
\]
where\vspace*{1pt} $S^2_X = \frac{1}{m-1}\sum_{i=1}^m (X_i - \bar{X}_{m})^2$ and
$S^2_Y = \frac{1}{n-1}\sum_{j=1}^n (Y_i - \bar{Y}_{m})^2$, can allow
one to construct a test that attains asymptotic rejection probability
$\alpha$ when $P \neq Q$ while providing an additional advantage of
maintaining \textit{exact} level $\alpha$ when $P = Q$.
\end{examplecontd}

\begin{examplecontd}[(Continued)]
Define the studentized median statistic
\[
S_{m,n} = \frac{N^{1/2}[\theta(\hat{P}_m)- \theta(\hat{Q}_n)]
}{\hat v_{m,n}},
\]
where $\hat v_{m,n}$ is a consistent estimator of $v(P,Q)$ defined
in (\ref{equationvvar}). There are several choices for a consistent
estimator of $v(P,Q)$. Examples include the usual kernel estimator
(Devroye and Wagner \cite{r4}), bootstrap estimator (Efron~\cite{r5}), and the smoothed bootstrap (Hall, DiCiccio, and Romano \cite{r6}).
\end{examplecontd}

\begin{remark}\label{remarklocalpower}
Suppose that the true unconditional distribution of a test $T_{m,n}$
is, under the null hypothesis, asymptotically given by a distribution
$R ( \cdot)$.
Typically a test rejects when $T_{m,n} > r_{m,n}$, where $r_{m,n}$
is nonrandom, as happens in many classical settings. Then, we typically
have $r_{m,n} \to r ( 1- \alpha) \equiv R^{-1} ( 1- \alpha)$.
Assume that $T_{m,n}$ converges to some limit law $R' ( \cdot)$ under
some sequence of alternatives which are contiguous to some distribution
satisfying the null. Then, the power of the test against such a sequence
would tend to $1- R' ( r(1- \alpha))$. The point here is that, under
the conditions of Theorem~\ref{theoremgeneral}, the permutation test
based on a random critical value $\hat r_{m,n}$ obtained from
the permutation\vspace*{-1pt} distribution satisfies, under the null, $\hat r_{m,n}
\Parrow r ( 1- \alpha)$.
But then, contiguity implies the same behavior under a sequence of contiguous
alternatives. Thus, the permutation test has the same limiting local power
as the ``classical'' test which uses the nonrandom critical value.
So, to first order, there is no loss in power in using a permutation
critical value.
Of course, there are big gains because the permutation test applies much
more broadly than for usual parametric models, in that it retains the
level exactly across a broad class of distributions and is at least
asymptotically
justified for a large nonparametric family.
\end{remark}

\section{Generalizations} \label{sectiongeneralization}

\subsection{Wilcoxon statistic and general $U$-statistics} \label{secustat}

So far, we considered two-sample problems where the statistic is based
on the difference of estimators that are asymptotically linear.
Although this class of estimators includes many interesting cases such
as testing equality of means, medians, and variances, it does not
include other important statistics like the Wilcoxon statistic or some
rank statistics where the parameter of interest is a function of the
joint distribution $\theta(P,Q)$ and not just a simple difference
$\theta(P) - \theta(Q)$.

In our companion paper (Chung and Romano \cite{r1}), however, we
consider these statistics in a more general $U$-statistic framework.
More specifically, assume that $X_1,\ldots, X_m$ are i.i.d. $P$, and
independently, $Y_1,\ldots, Y_n$ are i.i.d.~$Q$. The problem studied
is to test the null hypothesis
\[
H_0\dvtx  \E_{P,Q} \bigl( \varphi(X_1,\ldots,
X_r, Y_1,\ldots, Y_r) \bigr) = 0,
\]
which can be estimated by its corresponding two-sample $U$-statistic of
the form
\[
U_{m,n}(Z) = \frac{1}{{m \choose r}{n \choose r}} \sum_{\alpha}
\sum_{\beta} \varphi(X_{\alpha_1},\ldots,
X_{\alpha_r}, Y_{\beta_1},\ldots, Y_{\beta_r}),
\]
where $\alpha$ and $\beta$ range over the sets of all unordered
subsets of $r$ different elements chosen from $\{1,\ldots, m \}$ and
of $r$ different elements chosen from $\{1,\ldots, n \}$, respectively.

This general class of $U$-statistics covers, for example, Lehmann's
two-sample $U$-statistic to test $H_0\dvtx  P(|Y'-Y| > |X'-X|) = 1/2$, the
two-sample Wilcoxon statistic to test $H_0\dvtx  P( X \leq Y) = P(Y \leq
X)$, and some other interesting rank statistics. Under quite weak
assumptions, we provide a general theory whereby one can construct a
permutation test of a parameter $\theta(P, Q) = \theta_0$ which
controls the asymptotic probability of a type 1 error in large samples
while retaining the exactness property in finite samples when the
underlying distributions are identical. The technical arguments
involved in this $U$-statistic problem are different from Section~\ref
{secstudentized}, but the mathematics and statistical foundations to
be laid out in Section~\ref{sectioningredients} provide fundamental
ingredients that aid our asymptotic derivations.

\subsection{Robust $k$-sample test}\label{seconeway}

From our general considerations, we are now guided by the principle
that the large sample distribution of the test statistic should not
depend on the underlying distributions; that is, it should be asymptotically
pivotal under the null. Of course, it can be something other than normal,
and we next consider the important problem of testing equality of
parameters of $k$-samples (where a limiting Chi-squared distribution is
obtained).

Assume we observe $k$ independent samples of i.i.d. observations.
Specifically, assume $X_{i,1},\ldots, X_{i,n_i}$ are i.i.d.\vadjust{\goodbreak} $P_i$.
Some of our results will hold for fixed $n_1,\ldots, n_k$, but
we also have asymptotic results as $N \equiv\sum_i n_i \to\infty$.
Let $n = (n_1,\ldots, n_k )$, and the notation $n \to\infty$ will mean
$\min_i n_i \to\infty$. Let $\theta(\cdot)$ be a real-valued
parameter, defined on some space of distributions $\cal P$. The problem
of interest is to test the null hypothesis
%
\begin{equation}
\label{equationmultinull} H_0\dvtx  \theta(P_1) = \cdots=
\theta(P_k)
\end{equation}
against the alternative
\[
H_1\dvtx  \theta(P_i) \neq\theta(P_j)
\qquad\mbox{for some } i, j.
\]
When $P_1 = \cdots= P_k$ holds, one can construct permutation tests
with exact level $\alpha$. However, if $P_i \neq P_j$ for some $i, j$,
then the test may fail to achieve the rejection probability equal to
$\alpha$ even asymptotically.

We will assume that asymptotically linear estimators are available,
that is, (\ref{equationxlinear}) holds for i.i.d. samples under $P_i$
for $i = 1,\ldots, k$, where $f_{P_i}$ can depend on the underlying
distribution $P_i$. Further assume that the expansion also holds for
i.i.d. observations $\bar Z_{i,1},\ldots, \bar Z_{i, n_i}$ sampled
from the mixture distribution $\bar P = \sum_{i=1}^k p_i P_i$, where
$n_i/N \to p_i$. Note that the asymptotic linearity conditions need not
require any form of differentiability (though of course, some form of
differentiability is a sufficient condition). We will argue that the
asymptotic linearity conditions under $P_i$ for $i=1,\ldots, k$ and
$\bar P$, are sufficient to derive the asymptotic behavior of the
$k$-sample permutation distribution based on $T_{n,1}$ (defined below),
without having to impose any differentiability conditions.

The goal here is to construct a method that retains the exact control
of the probability of a type 1 error when the observations are i.i.d., but
also asymptotically controls the probability of a type 1 error under
very weak assumptions, specifically finite nonzero variances of the
influence functions.

\begin{lemma}\label{lemmamultichik}
Consider the above set-up. Assume (\ref{equationxlinear}) holds for
$P_1,\ldots, P_k$ with $0 < \sigma_i^2 = \sigma_i^2(f_{P_i}) = \E
_{P_i} f^2_{P_i}(X_{i,j}) < \infty$. Assume $n_i \to\infty$ with $n_i /
N \to p_i > 0$ for $i = 1,\ldots, k$. Let
%
\begin{equation}
\label{equationmultitn0} T_{n,0} = \sum
_{i=1}^k \frac{n_i}{\sigma_i^2} \biggl[ \hat
\theta_{n,i} - \frac{ \sum_{i=1}^k n_i \hat\theta_{n,i} /
\sigma_i^2 }{\sum_{i=1}^k n_i/ \sigma_i^2 } \biggr]^2,
\end{equation}
where $\hat\theta_{n,i} = \hat\theta_{n,i}(X_{i,1},\ldots, X_{i,
n_i})$ and $\sigma^2_i = \sigma_i^2(f_{P_i}) = \E_{P_i} f_{P_i}^2
(X_{i,j})$. Further assume that $\hat\sigma_{n,i} \equiv\hat\sigma
_{n,i}(X_{i,1},\ldots, X_{i, n_i})$ is a consistent estimator of
$\sigma_i = \sigma_i(f_{P_i})$ when $X_{i,1},\ldots, X_{i, n_i}$
are i.i.d. $P_i$, for $i =1,\ldots, k$. Define
%
\begin{equation}
\label{equationmultitn1} T_{n,1} = \sum
_{i=1}^k \frac{n_i}{\hat\sigma_{n,i}^2} \biggl[ \hat
\theta_{n,i} - \frac{ \sum_{i=1}^k n_i \hat\theta_{n,i} /
\hat\sigma_{n,i}^2 }{\sum_{i=1}^k n_i/ \hat\sigma_{n,i}^2 } \biggr]^2.
\end{equation}

Then, under $H_0$,
both $T_{n,0}$ and $T_{n,1}$ converge in distribution to the Chi-squared
distribution with $k-1$ degrees of freedom.
\end{lemma}

Let $\hat R_{n,1} ( \cdot)$ denote the permutation distribution
corresponding to $T_{n,1}$. In words, $T_{n,1}$ is recomputed over all
permutations of the data.
Specifically, if we let
\[
(Z_1,\ldots, Z_N ) = ( X_{1,1},\ldots,
X_{1, n_1}, X_{2,1},\ldots, X_{2, n_2},\ldots,
X_{k,1},\ldots, X_{k, n_k } ),
\]
then, $\hat R_{n,1} (t)$ is formally equal to the right-hand side of
(\ref{equationnewperm}), with $T_{m,n}$ replaced by $T_{n,1}$.

\begin{theorem}\label{theoremmultiperm}
Assume the same setup and conditions of Lemma~\ref{lemmamultichik}
with $0 < \sigma_i^2 = \sigma_i^2(f_{P_i}) = \E_{P_i}
f^2_{P_i}(X_{i,j}) < \infty$.
Assume $n_i \to\infty$ with $n_i / N \to p_i > 0$. Further\vspace*{1pt} assume
that the consistency of $\hat\sigma_{n,i}$ of $\sigma_i$ under $P_i$
also holds under $\bar P$ as well so that, when the $\bar Z_i$ are
i.i.d. $\bar P$,
\[
\hat\sigma_{n,i}(\bar Z_1,\ldots, \bar Z_{n_i})
\Parrow\sigma(f_{\bar P}) \qquad\mbox{as } n \to\infty
\]
with $0 < \sigma^2(f_{\bar P}) < \infty$.

Then, under $H_0$,
%
\begin{equation}
\label{equationmultiperm} \hat R_{n,1} (t) \Parrow G_{k-1}
(t),
\end{equation}
where $G_d$ denotes the Chi-squared distribution with $d$ degrees of freedom.
Moreover, if $P_1,\ldots, P_k$ satisfy $H_0$, then the probability
that the permutation test rejects $H_0$ tends to the nominal level
$\alpha$.
\end{theorem}

\begin{example}[(Nonparametric $k$-sample
Behrens--Fisher problem)]\label{exkbehrensfisher}
Consider the special case where $\theta_i (P) = \mu_i (P)$
is the population mean. Also, let $\hat\theta_{n,i}$ be the sample
mean of the $i$th sample.
When the populations are assumed normal with possibly different unknown
variances, this is
the classical Behrens--Fisher problem. Here, we do not assume normality
and provide a general solution
for testing the equality of parameters of several distributions.
Indeed, we have exact finite
sample type 1 error control when all the populations are the same, and
asymptotically type 1 error
control when the populations are possibly distinct. (Some relatively
recent large sample approaches
which do not retain our finite sample exactness property to this
specific problem are given in Rice and Gaines \cite{r25} and
Krishnamoorthy, Lu and Mathew \cite{r15}.)
\end{example}

\section{Simulation results}\label{secsimulation}
Monte Carlo simulation studies illustrating our results are presented
in this section. Table~\ref{tbmedians} tabulates the rejection
probabilities of one-sided tests for the studentized permutation median
test where the nominal level considered is $\alpha= 0.05$. The
simulation results confirm that the studentized permutation median test
is valid in the sense that it approximately attains level $\alpha$ in
large samples.

%
\begin{table}
\tabcolsep=0pt
\caption{Monte Carlo simulation results for studentized permutation
median test\break (one-sided, $\alpha= 0.05$)}\label{tbmedians}
\begin{tabular*}{\tablewidth}{@{\extracolsep{\fill}}lc c  c d{1.4} d{1.4} c c c@{}}
\hline
&
\multicolumn{1}{r}{$\bolds{m}$\textbf{:}} & \multicolumn{1}{c}{\textbf{5}}
& \multicolumn{1}{c}{\textbf{13}} & \multicolumn{1}{c}{\textbf{51}} & \multicolumn{1}{c}{\textbf{101}}
& \multicolumn{1}{c}{\textbf{101}}
& \multicolumn{1}{c}{\textbf{201}} & \multicolumn{1}{c@{}}{\textbf{401}} \\
\multicolumn{1}{@{}l}{\textbf{Distributions}} &
\multicolumn{1}{r}{$\bolds{n}$\textbf{:}} & \multicolumn{1}{c}{\textbf{5}}
& \multicolumn{1}{c}{\textbf{21}} & \multicolumn{1}{c}{\textbf{101}} & \multicolumn{1}{c}{\textbf{101}}
& \multicolumn{1}{c}{\textbf{201}} & \multicolumn{1}{c}{\textbf{201}} & \multicolumn{1}{c@{}}{\textbf{401}} \\
\hline
$N(0,1)$  &
Not studentized &  0.1079 & 0.1524 & 0.1324 & 0.2309 & 0.2266 & 0.2266 & 0.2249\\
$N(0,5)$ & Studentized & 0.0802 & 0.1458 & 0.095 & 0.0615 & 0.0517 & 0.0517 & 0.0531\\
[4pt]
$N(0,1)$ & Not studentized & 0.0646 & 0.1871 & 0.2411 & 0.1769 & 0.1849 & 0.1849 & 0.1853 \\
$T(5)$ & Studentized & 0.0707 & 0.1556 & 0.0904 & 0.0776 & 0.0661 & 0.0661 & 0.0611 \\
[4pt]
$\operatorname{Logistic}(0,1)$ & Not studentized & 0.0991
& 0.1413 & 0.1237 & 0.2258 & 0.2233 & 0.2233 & 0.2261\\
$U(-10, 10)$ & Studentized & 0.0771 & 0.1249 & 0.0923
& 0.0686 & 0.0574 & 0.0574 & 0.0574 \\
[4pt]
$\operatorname{Laplace}(\ln2, 1)$ & Not studentized & 0.0420 & 0.0462
& 0.0477 & 0.048 & 0.0493 & 0.0461 & 0.0501 \\
$\exp(1)$ & Studentized & 0.0386 & 0.0422 & 0.0444 & 0.0502
& 0.0485 & 0.0505 & 0.0531 \\
\hline
\end{tabular*}
%
\end{table}

In the simulation, odd numbers of sample sizes are selected in the
Monte Carlo simulation for simplicity.\vadjust{\goodbreak} We consider several pairs of
distinct sample distributions that share the same median as listed in
the first column of Table~\ref{tbmedians}.
For each situation, 10,000 simulations were performed. Within a given
simulation,
the permutation test was calculated by randomly sampling 999 permutations.
Note that neither the exactness properties nor the asymptotic properties
are changed at all (as long as the number of permutations sampled tends
to infinity). For a discussion on stochastic approximations to the permutation
distribution, see the end of Section 15.2.1 in Lehmann and Romano \cite{r19}
and Section 4 in Romano \cite{r26}.
As is well known, when the underlying distributions of two distinct
independent samples are not identical, the permutation median test is
not valid in the sense that the rejection probability is far from the
nominal level $\alpha= 0.05$. For example, although a logistic
distribution with location parameter 0 and scale parameter 1 and a
continuous uniform distribution with the support ranging from $-$10 to 10
have the same median of 0, the rejection probability for the sample
sizes examined is between 0.0991 and 0.2261 and moves further away from
the nominal level $\alpha= 0.05$ as sample sizes increase.

In contrast, the studentized permutation test results in rejection
probability that tends to the nominal level $\alpha$ asymptotically.
We apply\vspace*{1pt} the bootstrap method (Efron~\cite{r5}) to estimate the
variance for the median $\frac{1}{4 f^2_P(\theta)}$ in the simulation
given by
\[
m \sum_{l=1}^{m} \bigl[X_{(l)} -
\theta(\hat{P}_m) \bigr]^2 \cdot\P\bigl(\theta\bigl(
\hat{P}^*_m\bigr) = X_{(l)} \bigr),
\]
where for an odd number $m$,
\begin{eqnarray*}
\P\bigl(\theta\bigl(\hat{P}^*_m\bigr) = X_{(l)} \bigr) &=&
\P\biggl(\Binomial\biggl(m, \frac{l-1}{m} \biggr) \leq\frac{m-1}{2}
\biggr)\\
&&{} - \P\biggl(\Binomial\biggl(m, \frac{l}{m} \biggr) \leq
\frac{m-1}{2} \biggr).
\end{eqnarray*}
As noted earlier, there exist other choices such as the kernel
estimator and the smoothed bootstrap estimator. We emphasize, however,
that using the bootstrap to obtain an estimate of standard error does
\textit{not} destroy the exactness of permutation tests under identical
distributions.

\section{Four technical ingredients}\label{sectioningredients}

In this section, we discuss four separate ingredients, from which
the main results flow. These results are separated out so they
can easily be applied to other problems and so that the main technical
arguments are highlighted.
The first two apply more generally to randomization tests, not just
permutation tests, and are stated as such.

\subsection{Hoeffding's condition}\label{sectionhoeffding}

Suppose data $X^n$ has distribution $P_n$ in $\mathcal{X}_n$, and
$\mathbf{G}_n$ is a finite group of transformations $g$ of $\mathcal
{X}_n$ onto itself.
For a given statistic $T_n = T_n (X^n)$, let $\hat R_n^T ( \cdot)$ denote
the randomization distribution of $T_n$, defined by
%
\begin{equation}
\label{eqrand} \hat R_n^T (t) = \frac{1}{|G_n|}\sum
_{g \in G_n} I \bigl\{T_n\bigl(gX^n
\bigr) \leq t\bigr\},
\end{equation}
where $|G_n|$ denotes the cardinality of $G_n$.
Hoeffding \cite{r8} gave a sufficient condition to derive the limiting behavior
of $\hat R_n^T ( \cdot)$. This condition is verified repeatedly
in the proofs, but we add the result that the condition is also necessary.

\begin{theorem}\label{theoremhoeffding}
Let $G_n$ and $G_n'$ be independent and uniformly distributed
over~$\mathbf{G}_n$ (and independent of $X^n$). Suppose, under $P_n$,
%
\begin{equation}
\label{equationcondi} \bigl(T_n\bigl(G_nX^n
\bigr),T_n\bigl(G'_nX^n\bigr)
\bigr) \dto\bigl(T, T'\bigr),
\end{equation}
where $T$ and $T'$ are independent, each with common c.d.f.
$R^T(\cdot)$.
Then, for all continuity points $t$ of $R^T ( \cdot)$,
%
\begin{equation}
\label{equationhatRn} \hat R_n^T (t) \Parrow
R^T ( t ).
\end{equation}
Conversely, if (\ref{equationhatRn}) holds for some limiting c.d.f. $
R^T ( \cdot)$ whenever $t$ is a continuity point, then (\ref
{equationcondi}) holds.
\end{theorem}

The reason we think it is important to add the necessity part of the result
is that our methodology is somewhat different than that of other authors
mentioned in the \hyperref[secintroduction]{Introduction}, who take a more conditional approach
to proving limit theorems. After all, the permutation distribution is
indeed a distribution conditional on the observed set of observations
(without regard to ordering). However, the theorem shows that a sufficient
condition is obtained by verifying an unconditional weak convergence
property. Nevertheless, simple
arguments (see the supplementary appendix \cite{r2}) show the
condition is indeed necessary and so taking such
an approach is not fanciful.

\subsection{Slutsky's theorem for randomization distributions}
\label{sectionslutsky}

Consider the general setup of Section~\ref{sectionhoeffding}.
The result below describes Slutsky's theorem in the context
of randomization distributions. In this context, the randomization
distributions are random themselves, and therefore the usual Slutsky's
theorem does not quite apply. Because of its utility in the proofs of
our main results, we highlight the statement.
Given sequences of statistics $T_n$, $A_n$ and $B_n$, let $\hat R_n^{AT
+B} ( \cdot)$ denote the randomization
distribution corresponding to the statistic sequence $A_n T_n + B_n$;
that is,
replace $T_n$ in (\ref{eqrand}) by $A_n T_n + B_n$, so
%
\begin{equation}
\label{eqraa} \hat{R}_n^{AT +B} (t) \equiv\frac{1}{|G_n|}
\sum_{g \in G_n} \I\bigl\{A_n \bigl(g
X^n \bigr) T_n \bigl(gX^n\bigr) +
B_n \bigl( g X^n \bigr) \leq t\bigr\}.
\end{equation}

\begin{theorem}\label{theoremslutsky} Let $G_n$ and $G_n'$ be
independent and uniformly distributed over~$\mathbf{G}_n$ (and
independent of $X^n$).
Assume $T_n$ satisfies (\ref{equationcondi}).
Also, assume
%
\begin{equation}
\label{equationsluta} A_n \bigl(G_n X^n
\bigr) \Parrow a
\end{equation}
and
%
\begin{equation}
\label{equationslutb} B_n \bigl(G_n X^n
\bigr) \Parrow b
\end{equation}
for constants $a$ and $b$.
Let $R^{aT+b} ( \cdot)$ denote the distribution of $aT+b$, where $T$
is the limiting random variable assumed in (\ref{equationcondi}).
Then
\[
\hat R_n^{AT+B} (t) \Parrow R^{aT+b} (t),
\]
if the distribution $R^{aT+b} ( \cdot)$ of $aT+b$ is continuous at $t$.
[Of course, $R^{aT+b} (t) = R^T ( \frac{t-b}{a} )$ if $a \ne0$.]
\end{theorem}

\subsection{A coupling construction}\label{sectioncoupling}

Consider the general situation where $k$ samples are observed from
possibly different distributions. Specifically,
assume for $i=1,\ldots, k$ that $X_{i,1},\ldots, X_{i, n_i}$
is a sample of $n_i$ i.i.d. observations from $P_i$. All $N \equiv\sum
_i n_i$
observations are mutually independent.
Put all the observations together in one vector
\[
Z = ( X_{1,1},\ldots, X_{1, n_1}, X_{2,1},\ldots,
X_{2, n_2},\ldots, X_{k, 1},\ldots, X_{k, n_k } ).
\]

The basic intuition driving the results concerning the behavior of
the permutation distribution stems from the following.
Since the permutation distribution considers the empirical distribution
of a statistic evaluated at all permutations of the data, it clearly
does not depend on the ordering of the observations.
Let $n_i / N$ denote the proportion of observations
in the $i$th sample, and let $p_i = \lim_{n_i \to\infty} n_i/N \in
(0,1)$. Assume
that $n_i \to\infty$ in such a way that
%
\begin{equation}
\label{equationannoy} p_i - \frac{n_i}{N} = O\bigl(
N^{-1/2}\bigr).
\end{equation}
Then the behavior of the permutation distribution based on $Z$ should
behave approximately like the behavior of the permutation distribution
based on a sample of $N$ i.i.d. observations $\bar Z = ( \bar Z_1,\ldots, \bar Z_N )$ from the mixture distribution $\bar P \equiv p_1
P_1 + \cdots+ p_k P_k$.
Of course,\vspace*{1pt} we can think of the $N$ observations generated from $\bar P$
arising out of a two-stage process: for $ i = 1,\ldots, N$, first
draw an index $j$ at random with probability $p_j$; then, conditional
on the outcome being $j$, sample $\bar Z_i$ from $P_j$.
However, aside from the fact
that the ordering of the observations in $Z$ is clearly that of $n_1$
observations from $P_1$, following by $n_2$ observations from $P_2$,
etc., the original sampling scheme is still only approximately like
that of sampling from $\bar P$. For example, the number of
observations $\bar Z_i$ out of the $N$ which are from $P_1$ is binomial
with parameters $N$ and $p_1$ (and so has mean equal to $p_1 N \approx n_1$),
while the number of observations from $P_1$ in the original
sample $Z$ is exactly $n_1$.

Along the same lines, let $\pi= ( \pi(1),\ldots, \pi(N) )$ denote
a random permutation of $\{ 1,\ldots, N \}$.
Then, if we consider a random permutation of both $Z$
and $\bar Z$, then the number of observations in the first $n_1$
coordinates of $Z$ which were $X_1$'s has the hypergeometric distribution,
while the number of observations in the first $n_1$ coordinates of
$\bar Z$
which were $X_1$'s is still binomial.

We can\vspace*{1pt} make a more precise statement by constructing a certain coupling
of $Z$ and $\bar Z$. That is, except for ordering, we can construct
$\bar Z$ to include almost the same set of observations as in $Z$. The
simple idea goes as follows.
Given $Z$, we will construct observations $\bar Z_1,\ldots, \bar
Z_N$ via the two-stage process as above, using the observations drawn
to make up the $Z_i$
as much as possible. First, draw an index $j$ among $\{1,\ldots, k\}$
at random with probability $p_j;$ then, conditionally on the outcome
being $j$, set $\bar Z_1 = X_{j,1}$. Next, if the next index $i$ drawn
among $\{1,\ldots, k\}$ at random with probability $p_i$ is different
from $j$ from which $\bar Z_1$ was sampled, then $\bar Z_2 = X_{i,1};$
otherwise, if $i=j$ as in the first step, set $\bar Z_2 = X_{j,2}$. In
other words, we are going to continue to use the $Z_i$ to fill in the
observations $\bar Z_i$. However, after a certain point, we will get
stuck because we will have already exhausted all the $n_j$ observations
from the $j$th population
governed by $P_j$. If this happens and an index $j$ was drawn again,
then just sample a new observation $X_{j, n_j+1}$ from $P_j$.
Continue in this manner so that as many as possible of the original
$Z_i$ observations
are used in the construction of $\bar Z$. Now, we have both $\bar Z$
and $Z$.
At this point, $\bar Z$ and $Z$ have many of the same observations in common.
The number of observations which differ, say $D$, is the (random)
number of
added observations required to fill up $\bar Z$.
(Note that we are obviously using the word ``differ'' here to mean the
observations
are generated from different mechanisms, though in fact there may be a positive
probability that the observations still are equal if the underlying
distributions have atoms. Still, we count such
observations as differing.)

Moreover, we can reorder the observations in $\bar Z$ by a permutation
$\pi_0$
so that $Z_i$ and $\bar Z_{\pi_0 (i)}$ agree for all $i$ except for
some hopefully
small (random) number $D$. To do this,
recall that $Z$ has the observations in order, that is, the first $n_1$
observations arose from $P_1$ and the next set of $n_2$ observations
came from $P_2$, etc. Thus, to couple $Z$ and $\bar Z$, simply put all
the observations in $\bar Z$ which came from $P_1$ first up to $n_1$.
That is, if the number of observations in $\bar Z$ from $P_1$ is
greater than or equal to $n_1$, then $\bar Z_{\pi(i)}$ for $i = 1,\ldots, n_1$ are filled with the observations in $\bar Z$ which came
from $P_1$, and if the number was strictly greater than $n_1$, put them
aside for now. On the other hand, if the number of observations in
$\bar Z$ which came from $P_1$ is less than~$n_1$, fill up as many of
$\bar Z$ from $P_1$ as possible, and leave the rest of the slots among
the first $n_1$ spots blank for now. Next, move onto the observations
in $\bar Z$ which came from $P_2$ and repeat the above procedure for
$n_1+1,\ldots, n_1+n_2$ spots; that is, we start filling up the spots
from $n_1+1$ as many of $\bar Z$ which came from $P_2$ as possible up
to $n_2$ of them. After going though all the distributions $P_i$ from
which each of observations in $\bar Z$ came, one must then complete the
observations in $\bar Z_{\pi_0};$ simply ``fill up'' the empty spots
with the remaining observations that have been put aside. (At this
point, it does not matter where each of the remaining observations gets
inserted; but, to be concrete, fill the empty slots by inserting the
observations which came from the index $P_i$ in chronological order
from when constructed.)
This permuting of observations in $\bar Z$ corresponds to a permutation
$\pi_0$
and satisfies
$Z_i = \bar Z_{\pi_0 (i)}$ for indices $i$ except for $D$ of them.

For example, suppose there are $k=2$ populations. Suppose that $N_1$ of
the $\bar Z$ observations came from $P_1$ and so
$N - N_1$ from $P_2$.
Of course, $N_1$ is random and has the binomial distribution with
parameters $N$
and $p_1$. If $ N_1 \ge n_1$, then the above construction yields
the first $n_1$ observations in $Z$ and $\bar Z_{\pi_0}$\vspace*{1pt} completely agree.
Furthermore, if $N_1 > n_1$, then the number of observations in $\bar
Z$ from $P_2$
is $N - N_1 < N - n_1 = n_2$, and $N-N_1$ of the last $n_2$ indices in $Z$
match those of $\bar Z_{\pi_0}$, with the remaining differ.
In this situation, we have
\[
Z = ( X_ 1,\ldots, X_{n_1}, Y_1,\ldots,
Y_{n_2} )
\]
and
\[
\bar Z_{\pi_0} = (X_1,\ldots, X_{n_1},
Y_1,\ldots, Y_{N - N_1}, X_{n_1+1},\ldots,
X_{N_1} ),
\]
so that $Z$ and $\bar Z_{\pi_0}$ differ only in the last $N_1 -n_1$ places.
In the opposite situation where $N_1 < n_1$, $Z$ and $\bar Z_{\pi}$
are equal in the first $N_1$ and last $n_2$ places, only differing in spots
$N_1+1,\ldots, n_1$.

The number of observations $D$ where $Z$ and $\bar Z_{\pi_0}$ differ
is random
and it can be shown that
%
\begin{equation}
\label{equationsized} E(D/N) \leq N^{-1/2};
\end{equation}
see supplementary appendix \cite{r2}. In summary, the coupling
construction shows that only a fraction of the $N$ observations in $Z$
and $\bar Z_{\pi_0}$ differ with high probability.
Therefore, if the randomization distribution is based on a statistic
$T_N ( Z )$ such that the difference $T_ N (Z) - T_N ( \bar Z_{\pi_0}
)$ is small in some sense
whenever $Z$ and $\bar Z_{\pi_0}$ mostly agree, then one should be able
to deduce the behavior of the permutation distribution under
samples from $P_1,\ldots, P_k$ from the behavior of the permutation
distribution
when all $N$ observations come from the same distribution $\bar P$.
Whether or not this can be done requires some knowledge of the form
of the statistic, but intuitively it should hold if the statistic
cannot strongly be
affected
by a change in a small proportion of the observations; its validity
though must
be established on a case by case basis. Although the assessment of the
validity needs to be taken on a case by case basis, it readily extends
to a broader class of statistics such as ``mean-like'' statistics.
(However, this coupling argument and the contiguity results in Section
\ref{sectioncontiguity} together allow us to prove quite general
results.) The point is that it is a worthwhile and beneficial route to
pursue because
the behavior of the permutation distribution under $N$ i.i.d.
observations is
typically much easier to analyze than under the more general setting
when observations have possibly different distributions.
Furthermore, the behavior under i.i.d. observations seems fundamental
as this is the requirement for the ``randomization hypothesis'' to hold,
that is,
the requirement to yield exact finite sample inference.

To be more specific, suppose $\pi$ and $\pi'$ are independent random
permutations, and independent of the $Z_i$ and $\bar Z_i$.
Suppose we can show that
%
\begin{equation}
\label{equationnadal} \bigl( T_N ( \bar Z_{\pi} ),
T_N ( \bar Z_{\pi' } )\bigr) \dto\bigl(T, T'
\bigr),
\end{equation}
where $T$ and $T$ are independent with common c.d.f. $R ( \cdot)$.
Then, by Theorem~\ref{theoremhoeffding},
the randomization distribution based on $T_N$ converges in probability
to $R ( \cdot)$ when all observations are i.i.d. according to $\bar P$.
But since $\pi\pi_0$ (meaning $\pi$ composed with $\pi_0$ so $\pi
_0$ is applied first) and $\pi' \pi_0$ are also
independent random permutations, (\ref{equationnadal}) also implies
\[
\bigl( T_N ( \bar Z_{\pi\pi_0 } ), T_N ( \bar
Z_{\pi' \pi_0 } )\bigr) \dto\bigl(T, T'\bigr).
\]
Using the coupling construction to construct $Z$, suppose it can be shown
that
%
\begin{equation}
\label{equationshowg} T_N ( \bar Z_{\pi\pi_0} ) -
T_N ( Z_{\pi} ) \Parrow0.
\end{equation}
Then, it also follows that
\[
T_N ( \bar Z_{\pi' \pi_0} ) - T_N ( Z_{\pi'} )
\Parrow0,
\]
and so by Slutsky's theorem,
it follows that
%
\begin{equation}
\label{equationcouple} \bigl( T_N ( Z_{\pi} ),
T_N ( Z_{\pi' } )\bigr) \dto\bigl(T, T'\bigr).
\end{equation}
Therefore, again by Theorem~\ref{theoremhoeffding}, the randomization
distribution also converges in probability to $R( \cdot)$
under the original model of $k$ samples from possibly different distributions.
In summary, the coupling construction of $\bar Z$, $Z$ and $\pi_0$
and the one added requirement (\ref{equationshowg}) allow us to reduce
the study of the permutation distribution under possibly $k$ different
distributions
to the i.i.d. case when all $N$ observations are i.i.d. according to
$\bar P$.
We summarize this as follows.

\begin{lemma}\label{lemmacoupling}
Assume (\ref{equationnadal}) and (\ref{equationshowg}).
Then (\ref{equationcouple}) holds, and so the permutation distribution
based on $k$ samples from possibly different distributions behaves
asymptotically as if all observations are i.i.d. from the mixture
distribution $\bar P$ and satisfies
\[
\hat R^T_{m,n} (t) \Parrow R(t),
\]
if $t$ is a continuity point of the distribution $R$ of $T$ in (\ref
{equationnadal}).
\end{lemma}

\begin{example}[(Difference of sample means)]\label{examplediffmeans}
To appreciate what is involved in the verification of (\ref{equationshowg}),
consider the two-sample problem considered in Theorem \ref
{theorembasic}, in the
special case of testing equality of means. The unknown variances
may differ and are assumed finite.
Consider
the test statistic $T_{m,n} = N^{1/2} [ \bar X_m - \bar Y_n ]$.
By the coupling construction, $\bar Z_{\pi\pi_0}$ and $Z_{\pi}$
have the same components except for at most $D$ places. Now,
\begin{eqnarray*}
T_{m,n} ( \bar Z_{ \pi\pi_0} ) - T_{m,n} ( Z_{\pi}
) &=& N^{1/2} \Biggl[ \frac{1}{m} \sum
_{i=1}^m ( \bar Z_{ \pi\pi_0 (i)} - Z_{\pi(i) }
)\Biggr] \\
&&{}- N^{1/2} \Biggl[ \frac{1}{n} \sum
_{j=m+1}^N ( \bar Z_{ \pi\pi_0 (j)} -
Z_{ \pi(j) }) \Biggr].
\end{eqnarray*}
All of the terms in the above two sums are zero except for at most $D$
of them.
But any nonzero term like $\bar Z_{\pi\pi_0 (i)} - Z_{\pi(i) }$ has variance
bounded above by
\[
2 \max\bigl( \operatorname{Var} (X_1 ), \operatorname{Var}(Y_1 ) \bigr) < \infty.
\]
Note the above random variable has mean zero under the null hypothesis
that $E(X_i ) = E( Y_j )$. To bound its variance, condition on\vadjust{\goodbreak} $D$ and
$\pi$,
and note it has conditional mean 0 and conditional variance
bounded above by
\[
N \frac{1}{\min(m^2, n^2 )} 2 \max\bigl(\operatorname{Var} (X_1 ), \operatorname{Var}
(Y_1)\bigr) D
\]
and hence unconditional variance bounded above by
\[
N \frac{1}{\min(m^2, n^2 )} 2 \max\bigl(\operatorname{Var} (X_1 ), \operatorname{Var}
(Y_1)\bigr) O \bigl( N^{1/2} \bigr) = O \bigl(
N^{-1/2} \bigr) = o(1),
\]
implying (\ref{equationshowg}). In words, we have shown that the behavior
of the permutation distribution can be deduced from the behavior of
the permutation\vspace*{1pt} distribution when all observations are i.i.d. with
mixture distribution $\bar P$.
\end{example}

Two final points are relevant. First, the limiting distribution $R$ is
typically the same as the limiting distribution of the true unconditional
distribution of $T_N$ under~$\bar P$. This is intuitively the case
because the permutation distribution is invariant under any permutation
of the combined data, and so the set of $N$ observations with exactly
$n_i$ observations sampled from $P_i$ and then randomly permuting them
behaves very nearly the same as a sample of $N$ observations from $\bar
P$. On the other hand, the true limiting distribution of the test
statistic under
$(P_1,\ldots, P_k )$ need not be the same as under $\bar P$ as it
will in general depend on the underlying distributions $P_1,\ldots, P_k$.
However, suppose the choice of test statistic $T_N$ is such that
it is an asymptotic pivot in the sense that its limiting distribution does
not depend on the underlying probability distributions.
Then, the limiting distribution of the test statistic will be the same
whether sampling from $(P_1,\ldots, P_k)$ or $(\bar P,\ldots, \bar
P)$. In such cases, the randomization or permutation distribution under
$(P_1,\ldots, P_k )$ will asymptotically
reflect the true unconditional distribution of $T_N$, resulting in
asymptotically
valid inference. Indeed, the general results in Section \ref
{secstudentized} yield many
examples of this phenomenon. However, that these statements need
qualification is made clear by the following two (somewhat contrived) examples.

\begin{example}\label{exampleex1}
Here, we illustrate a situation where coupling works, but the true
sampling distribution does not
behave like the permutation distribution under the mixture model $\bar P$.
In the two-sample setup with $m=n$, suppose $X_1,\ldots, X_n$ are
i.i.d.
according to uniformity on the set of $x$ where $|x| <1$, and $Y_1,\ldots, Y_n$
are i.i.d. uniform on the set of $y$ with $2 < |y| < 3$. So, $E(X_i ) =
E(Y_j ) = 0$.
Consider a test statistic $T_{n,n}$ defined as
\[
T_{n,n} ( X_1,\ldots, X_n, Y_1,\ldots, Y_n ) = N^{-1/2} \Biggl[ \sum
_{i=1}^n I \bigl\{ |Y_i | > 2 \bigr\} - I \bigl\{ | X_i | < 2 \bigr\} \Biggr].
\]
Under the true sampling scheme, $T_{n,n}$ is
zero with probability one. However, if all $2n$ observations are
sampled from
the mixture model, it is easy to see that $T_{n,n}$ is asymptotically normal
$N(0,1/4)$, which is the same limit for the permutation distribution
(in probability).
So here, the permutation distribution under the given distributions
is the same as under $\bar P$, though it does not reflect the actual true
unconditional sampling distribution.
\end{example}

\begin{example}\label{exampleex2}
Here, we consider a situation where both populations are indeed
identical, so there is no need for a coupling argument. However, the
point is that the permutation distribution does not behave like the
true unconditional sampling distribution. Assume $X_1,\ldots, X_n$ and
$Y_1,\ldots, Y_n$ are all i.i.d. $N(0,1)$ and consider the test
statistic
\[
T_{n,n} (X_1,\ldots, X_n, Y_1,\ldots, Y_n ) = N^{-1/2} \sum_{i=1}^n
(X_i + Y_i ).
\]
Unconditionally, $T_{n,n}$ converges in distribution to $N(0,1)$.
However, the permutation distribution places mass one at $\frac
{n}{\sqrt{N}} ( \bar X_n + \bar Y_n )$ because the statistic $T_{n,n}$
is permutation invariant.
\end{example}

Examples~\ref{exampleex1} and~\ref{exampleex2} show that the
intuition provided in the paragraph before Example~\ref{exampleex1}
does not always work. However, in the two examples, the test statistic
does not reflect an actual comparison between $P$ and $Q$. Of course,
our theorems apply to tests of equality of parameters, and therefore
the test statistics are based on appropriate differences.

\subsection{An auxiliary contiguity result}\label{sectioncontiguity}

Consider the general situation involving $k$ (possibly distinct)
populations\vspace*{1pt} for $i = 1,\ldots, k$ with $n_i$ observations from
population $i$. Set $N = \sum_{i=1}^k n_i$ and $n = (n_1,\ldots,
n_k)'$, where the notation $n \to\infty$ means $\min_i n_i \to
\infty$. Assume all $N$ observations are mutually independent. Define
$p_{n,i} = n_i/N \to p_i \in(0,1)$ as $n_i \to\infty$ for $i = 1,\ldots, k$. Let $P_n$ be the multinomial distribution based on
parameters $s = s(n)$ and $p_n = (p_{n,1},\ldots, p_{n,k})$. So,
under $P_n$, let $M_{n,j}$ be the number of observations of type $i$
when $s$ observations are taken with replacement from a population with
$n_i$ observations of type $i$. So, $M_n \equiv(M_{n,1},\ldots,
M_{n,k}) \sim P_n$. Also, let $Q_n$ be the multivariate hypergeometric
distribution. Under $Q_n$, let $H_{n,i}$ be the number of observations
of type $i$ when $s$ observations are taken without replacement. So,
$H_n \equiv(H_{n,1},\ldots, H_{n,k}) \sim Q_n$.

We shall show that the multinomial distribution $P_n$ and the
multivariate hypergeometric distribution $Q_m$ are mutually contiguous,
which will allow us to obtain the limiting behavior of a statistic
under the given samples from $k$ probability distributions $P_i$ for $i
= 1,\ldots, k$, by instead calculating the limiting behavior of the
statistic when all $N$ observations are i.i.d. from the mixture
distribution $\bar P = \sum_{i=1}^k p_iP_i$, which is relatively
easier to obtain. For basic details on contiguity, see Section 12.3 in~Lehmann and Romano~\cite{r19}.

\begin{lemma}\label{lemmamulticontig}
Assume the above setup with $s/N \to\theta\in[0,1)$ as $n \to\infty$.
Consider the likelihood ratio $L_n (x) = dQ_n (x)/
dP_n (x) $.

\begin{longlist}
\item
The limiting distribution of $L_n (M_n)$ satisfies
%
\begin{equation}
\label{equationmultibasic} L_n (M_n ) \Larrow(1-\theta
)^{-({k-1})/{2}} \exp\biggl\{ - \frac{\theta}{2(1-\theta)} \chi^2_{k-1}
\biggr\},
\end{equation}
where $\chi^2_{k-1}$ denotes the Chi-squared distribution with $k-1$
degrees of freedom.

\item $Q_n$ and $P_n$ are mutually contiguous.
\end{longlist}
\end{lemma}

\begin{remark}\label{remarkmultinomial}
With $M_n \equiv(M_{n,1},\ldots, M_{n,k})$ having the multinomial
distribution with parameters $s$ and $p_n = (p_{n,1},\ldots, p_{n,k})$
as in Lemma~\ref{lemmamulticontig}, also let $\bar M_n \equiv(\bar
M_{n,1},\ldots, \bar M_{n,k})$ have the multinomial
distribution with parameters $s$ and $p = (p_1,\ldots, p_k)$. Then,
the distributions
of $M_n$ and $\bar M_n$ are contiguous if and only if
$p_{n,i} - p_i = O ( n_i^{-1/2} )$, not just $p_{n,i} \to p_i$ for all
$i= 1,\ldots, k$.
\end{remark}

\begin{lemma}\label{lemmamultibasicp}
Suppose $V_1,\ldots, V_s$ are i.i.d. according to the mixture distribution
\[
\bar P \equiv\sum_{i=1}^k
p_iP_i,
\]
where $p_i \in(0,1), \sum_{i=1}^k p_i = 1$ and $P_i$'s are probability
distributions on some general space. Assume, for some sequence $W_n$ of
statistics,
%
\begin{equation}
\label{equationmultiTmV} W_n ( V_1,\ldots,
V_s ) \Parrow t
\end{equation}
for some constant $t$ (which can depend on the $P_i$'s and $p_i$'s).
Let $n_i \to\infty$, \mbox{$s(n) \to\infty$}, with $s/N \to\theta\in[0, 1)$,
$N = \sum_{i=1}^k n_i$, $p_{n,i} = n_i/N$, and $p_{n,i} \to p_i \in
(0,1)$ with
%
\begin{equation}
\label{equationmultimn} p_{n,i} - p_i = O\bigl(
n_i^{-1/2} \bigr).
\end{equation}
Further, let $X_{i,1},\ldots, X_{i,n_i}$ be i.i.d. $P_i$ for $i = 1,\ldots, k$.
Let
\[
(Z_1,\ldots, Z_N) = ( X_{1,1},\ldots,
X_{1,n_1},\ldots, X_{k,
1},\ldots, X_{k, n_k} ).
\]
Let $( \pi(1),\ldots, \pi(N) )$ denote a random permutation
of $\{ 1,\ldots, N \}$ (and independent of all other variables).
Then,
%
\begin{equation}
\label{equationmultiTmZ} W_n ( Z_{\pi(1)},\ldots,
Z_{\pi( s)} ) \Parrow t.
\end{equation}
\end{lemma}

\begin{remark}
The importance of Lemma~\ref{lemmamultibasicp} is that is allows us to
deduce the behavior of the statistic $W_n$ under the randomization or
permutation distribution from the basic assumption\vadjust{\goodbreak} of how $W_n$ behaves
under i.i.d. observations from the mixture distribution $\bar P$. Note
that in (\ref{equationmultiTmV}), the convergence in probability
assumption is required when the $V_i$ are $\bar P$ (so the $P$ over the
arrow is just a generic symbol for convergence in probability).
\end{remark}

\section{Conclusion}

When the fundamental assumption of identical distributions need not
hold, two-sample permutation tests are invalid unless quite stringent
conditions are satisfied depending on the precise nature of the
problem. For example, the two-sample permutation test based on the
difference of sample means is asymptotically valid only when either the
distributions have the same variance or they are comparable in sample
size. Thus, a careful interpretation of rejecting the null is
necessary; rejecting the null based on the permutation tests does not
necessarily imply a valid rejection of the null that some real-valued
parameter $\theta(F,G)$ is some specified value $\theta_0$. We
provide a framework that allows one to obtain asymptotic rejection
probability $\alpha$ in two-sample permutation tests. One great
advantage of utilizing the proposed test is that it retains the
exactness property in finite samples when $P = Q$, a desirable property
that bootstrap and subsampling methods fail to possess.

To summarize, if the true goal is to test whether the parameter of
interest $\theta$ is some specified value $\theta_0$, permutation
tests based on correctly studentized statistic is an attractive choice.
When testing the equality of means, for example, the permutation
$t$-test based on a studentized statistic obtains asymptotic rejection
probability $\alpha$ in general while attaining exact rejection
probability equal to $\alpha$ when $P = Q$. In the case of testing the
equality of medians, the studentized permutation median test yields the
same desirable property. Moreover, the results extend to quite general
settings based on asymptotically linear estimators.
The results extend to $k$-sample problems as well, and analogous results
hold in $k$-sample problem of comparing general parameters, which
includes the nonparametric $k$-sample Behrens--Fisher problem.
The guiding principle is to use a test statistic that is asymptotically
distribution-free
or pivotal. Then, the technical arguments developed in this paper can be
shown that the permutation test behaves asymptotically the same as when
all observations share a common distribution. Consequently, if the permutation
distribution reflects the true underlying sampling distribution, asymptotic
justification is achieved.

As mentioned in the \hyperref[secintroduction]{Introduction}, proper implementation of a
permutation test is vital
if one cares about confirmatory inference through hypothesis testing;
indeed, proper
error control of types 1, 2 and 3 errors can be obtained for test of parameters
by basing inference on test statistics which are asymptotically pivotal.
Thus, the foundations are laid for considering more complex problems in
modern data
analysis, such as two-sample microarray genomics problems, where a very large
number of tests are performed simultaneously. (Indeed, there are many
microarray analyses
which have begun by performing a permutation test for each gene,
without proper studentization.)
The role of permutations in multiple testing cannot be properly
understood without
a firm basis for single testing.
Thus, future work will further develop the ideas presented here so that
permutation tests can be
applied to other measures of error control in multiple testing such as
the false discovery rate.


\begin{supplement}
\stitle{Supplement to ``Exact and asymptotically robust permutation tests''}
\slink[doi]{10.1214/13-AOS1090SUPP} 
\sdatatype{.pdf}
\sfilename{aos1090\_supp.pdf}
\sdescription{Contains
proofs of all the results in the paper.}
\end{supplement}


\printaddresses

\end{document}